\documentclass[letterpaper, 10 pt, conference]{ieeeconf} 
\IEEEoverridecommandlockouts      
\usepackage[T1]{fontenc}
\usepackage{mathptmx}
\usepackage{amsmath,amssymb,mathtools,graphicx,booktabs,cite,balance}

\usepackage{amsthm}
\usepackage[hidelinks]{hyperref}
\hypersetup{pdftitle={Quadratic Returns for a Moving Water Tank: An Explicit Right Inverse and a Sylvester Representation},pdfauthor={Constantinos Kitsos}}

\theoremstyle{plain}
\newtheorem{theorem}{Theorem}
\newtheorem{proposition}[theorem]{Proposition}
\newtheorem{lemma}[theorem]{Lemma}
\newtheorem{corollary}[theorem]{Corollary}

\theoremstyle{definition}

\newtheorem{remark}[theorem]{Remark}

\newcommand{\E}{\mathcal E}
\newcommand{\Q}{\mathcal Q}
\newcommand{\Rq}{\mathcal R_{\rm q}}

\newcommand{\reflec}{\mathsf R}

\title{\LARGE \bf
Quadratic Returns for a Moving Water Tank: An Explicit Right Inverse and a Sylvester Representation
}

\author{Constantinos Kitsos
}

\begin{document}
\maketitle
\thispagestyle{empty}
\pagestyle{empty}
\begin{abstract}
Coron's return method established local nonlinear controllability of a moving water tank in sufficiently large time. The quadratic response at uniform equilibrium is studied in this work, where the acceleration input in the linearized model cannot excite the Fourier modes of the fluid with even indices. At three one-way wave crossings ($T=3$), the quadratic final-state map on a family of delayed controls that return the linearized state to zero is reduced to a scalar expression associated with a planar curve. Making this scalar expression constant provides a nonzero quadratic return, and a tangential variation yields an explicit right inverse of the derivative of the quadratic final-state map. This inverse acts on all missing modes and maps mean-zero periodic $H^r$ targets to $H^r$ controls for every nonnegative integer $r$. A fixed-point iteration also solves the quadratic final-state equation for sufficiently small targets. Finally, a time-dependent Sylvester equation provides an alternative representation of the coupling responsible for the second-order correction. The associated change of state coordinates involves one spatial derivative, while the control correction preserves regularity in appropriate Sobolev spaces. These results concern the first two variation equations and do not establish nonlinear controllability at $T=3$.
\end{abstract}

 \section{Introduction}
A tank moved by a single horizontal actuator can develop various fluid oscillations. The problem is relevant to the transport of partially filled containers, where the motion of the container also excites the liquid. For the inviscid Saint-Venant model at a uniform equilibrium, considered here, the  acceleration (control input) in the linearized model excites only the spatial Fourier modes with odd indices. The modes with even indices can, therefore, be controlled only through nonlinear interactions.

In \cite{DuboisPetitRouchon1999,PetitRouchon2002}, explicit motion planning for the linearized tank was developed. In \cite[Thm. 2]{Coron}, Coron established local controllability of the nonlinear shallow-water model (inviscid Saint-Venant equations) in sufficiently large time by the return method. An obstruction below two one-way wave crossings ($T<2$) was later proved in \cite{CKN}, whose unfolding and quadratic analysis are used here. However, this lower bound does not imply controllability at every longer time. Related stabilization results concern the inviscid tank linearized around nonuniform states \cite{CHXZ} and a nonlinear viscous model \cite{KarafyllisKrstic2022}, but these concern different models or reference trajectories from the quadratic variation equations studied here.

Inspired by the delayed control family constructed in \cite{PetitRouchon2002}, we derive a scalar formulac for the quadratic final-state map at $T=3$. We construct a nonzero smooth quadratic return and an explicit right inverse of the derivative of this map at that return, with $H^r$-to-$H^r$ bounds on all missing modes for every nonnegative integer $r$. A convergent fixed-point iteration then solves the quadratic final-state equation for small targets. Related quadratic recovery mechanisms appear in other problems for partial differential equations (PDEs) \cite{Cerpa2007,BeauchardCoron2006,BeauchardMarbach2020}. The time $T=3$ accommodates two successive pieces of an auxiliary function and a unit delay (no optimality for this time is claimed here).

A complementary time-dependent Sylvester transformation describes the same coupling in the variation equations. The differential transformation, subject to Sylvester operator equations, in \cite[Sec. 3.3.2, Eq. (3.88)]{KitsosThesis}, is associated with the hyperbolic observer design of \cite{KBP}, where spatial derivatives of the measured component enter the construction. Relevant differential transformations are used in \cite{KitsosBajodekBaudouin}, while the modal approach in \cite{KitsosFridman2024} avoids spatial differentiation in the implemented feedback for underactuated reaction--diffusion systems in cascade.  Sylvester equations are used in observer design, ODE--PDE compensation, and forwarding \cite{Emirsajlow,Natarajan,MarxBrivadisAstolfi}, with connections to backstepping kernels
\cite{HumalojaBekiarisLiberis2026}. For the present quadratic cascade, this representation
distinguishes the derivative requirement of the state transformation from the control correction, which preserves regularity.

Section \ref{sec:model} states the main result. Section \ref{sec:curve} presents the scalar reduction, reference curve, and right inverse. Section \ref{sec:sylvester} derives the Sylvester representation, and Section \ref{sec:numerics} illustrates a numerical example.

\noindent\textit{Notation.}
Sobolev space $H^r(I)$ is equipped with its usual Sobolev norm for an interval $I$ and a nonnegative integer $r$, and $C_c^\infty(I)$ denotes smooth functions with compact support in $I$. Primes denote derivatives with respect to the displayed argument, dots denote time derivatives, and subscripts $t,x$ denote partial derivatives.

\section{Model and main result}\label{sec:model}
\subsection{The first two variations}\label{subsec:firstvariations}
The inviscid Saint-Venant model on $[0,T]\times [0,1]$
\begin{equation}\label{eq:model}
\begin{aligned}
H_t+(Hv)_x&=0,&v_t+H_x+vv_x&=-u(t),\\ v(t,0)&=v(t,1)=0,&\dot s&=u,\quad\dot d=s
\end{aligned}
\end{equation}
is considered here, with $H$ the water depth, $v$ the water velocity relative to the tank, $s$ and $d$ the tank velocity and position, and $u$ the tank acceleration corresponding to the control input. The state $(H,v,s,d)$ starts from initial condition $(1,0,0,0)$ and satisfies, by virtue of mass conservation, $\int_0^1H(t,x) dx=1$. Note that, for simplicity, since we use dimensionless variables, one time unit corresponds to the physical wave crossing time $L/\sqrt{g_{\rm grav}H_{\rm eq}}$, where $L$ is tank length, $g_{\rm grav}$ is gravitational acceleration, and $H_{\rm eq}$ is equilibrium depth.

Injecting in \eqref{eq:model}  an input of the form $u=\alpha f$, we may use the formal expansion
\begin{equation}\label{eq:expansion}
H=1+\alpha h_1+\alpha^2h_2+\cdots,\qquad v=\alpha v_1+\alpha^2v_2+\cdots,
\end{equation}
leading to the variation equations
\begin{align}
h_{1,t}+v_{1,x}&=0,&v_{1,t}+h_{1,x}&=-f,\label{eq:first}\\
h_{2,t}+v_{2,x}&=-(h_1v_1)_x,& v_{2,t}+h_{2,x}&=-(v_1^2/2)_x,\label{eq:second}
\end{align}
as in \cite[Sec. 1.2]{CKN}. Both pairs have zero initial data and satisfy $v_i(t,0)=v_i(t,1)=0$, for smooth compatible inputs. These equations are studied directly, avoiding, in the present context, any nonlinear remainder estimate. The tank variables are exactly linear in the input $\alpha f$, satisfying $s=\alpha s_1$, $d=\alpha d_1$ exactly, where $\dot s_1=f$, $\dot d_1=s_1$, with zero initial data. If one takes $L^2$ inputs, the variation equations are understood through the unfolded transport formulation below, which incorporates the wall conditions by reflection and periodic continuation.

To identify the modes that do not appear in the response of order one, one may expand the height perturbation and the fluid velocity in $L^2(0,1)$ bases as $h_1(t,x)=\sqrt2\sum_{n\ge1}a_n(t)\cos(n\pi x),$ $v_1(t,x)=\sqrt2\sum_{n\ge1}b_n(t)\sin(n\pi x)$ and, then, linear equations \eqref{eq:first} yield
\begin{equation}\label{eq:modes}
\begin{gathered}
\dot a_n=-n\pi b_n,\qquad \dot b_n=n\pi a_n-\beta_n f;\\
\beta_n=\frac{\sqrt2(1-(-1)^n)}{n\pi}.
\end{gathered}
\end{equation}
We observe that since $\beta_n=0,$ for any $n$ even, all modes with even indices are not affected by the linear input, whereas the quadratic terms in \eqref{eq:second} provide the couplings through which these modes can be reached. 

Let us fix final time $T=3$ and define final-state maps by
\begin{equation}\label{eq:final-coefficients}
\begin{aligned}
\E_1[f]&=(h_1(T,\cdot),v_1(T,\cdot),s_1(T),d_1(T)),\\
\E_2[f]&=(h_2(T,\cdot),v_2(T,\cdot),0,0).
\end{aligned}
\end{equation}
A \emph{linear return} satisfies $\E_1[f]=0$, while a \emph{quadratic return} satisfies both $\E_1[f]=0$ and $\E_2[f]=0$. These conditions concern the final-state coefficients, without requiring that the fluid remains at equilibrium during the motion.

Following \cite[Def. 2.1 and Prop. 2.2]{CKN}, let us unfold each fluid pair as
\begin{equation}\label{eq:unfold}
z_i(t,x)=
\begin{cases}
h_i(t,x)+v_i(t,x),&0<x<1,\\
h_i(t,-x)-v_i(t,-x),&-1<x<0,
\end{cases}
\end{equation}
and then extend it with period $2$ in $x$. It turns out that $z_1(t,x+1)=-z_1(t,x)$, whereas $z_2$ is $1$-periodic and has zero mean \cite[Rems. 2.4 and 3.2]{CKN}. We also write $z_i[f]$ for the unfolded response of order $i$ corresponding to the input $f$. If we define quadratic final-state  map by
\begin{equation}\label{eq:Q}
\Q(f):=z_2[f](T,\cdot)|_{(0,1)},
\end{equation}
it follows that $\Q(f)=0$ is equivalent to $\E_2[f]=0$. When evaluated outside $(0,1)$, $\Q(f)$ denotes its $1$-periodic extension. If $\Q(f)=q$, the corresponding second-order height and velocity at final time are 
\begin{equation}\label{eq:recovered}
h_2(T,x)=\frac{q(x)+q(1-x)}2, v_2(T,x)=\frac{q(x)-q(1-x)}2. \end{equation}

Let now $\mathcal H^r(0,1)$, with $r$ a nonnegative integer, consist of $q\in H^r(0,1)$ with zero mean and matching traces $q^{(j)}(0)=q^{(j)}(1)$ for $0\le j<r$, and we equip this space with the usual $H^r(0,1)$ norm. We also define $\mathcal H^0(0,1)=L^2_0(0,1):=\{q\in L^2(0,1):\int_0^1q(x) dx=0\}$. In general, targets $q$ are extended with period $1$ whenever evaluated outside $(0,1)$.

\subsection{A family of linear returns}
Let $G\in H^2_0(0,2)$ be an auxiliary function of time, where $H^2_0(0,2)$ consists of $G\in H^2(0,2)$ satisfying $G(0)=G'(0)=G(2)=G'(2)=0$. Function $G$ is extended by zero outside $(0,2)$ and we define the input (corresponding to tank acceleration) by
\begin{equation}\label{eq:f_G}
f_G(t):=G''(t)+G''(t-1),\qquad 0<t<3.
\end{equation}
This is a shifted form of the delayed parametrization in \cite[Sec. II-A, Lemma 3]{PetitRouchon2002}, see also \cite[Lemma 2.7, proof]{CKN}, which credits \cite[Sec. 3.4]{DuboisPetitRouchon1999}. The trajectory below shows that every such control $f_G$ is a linear return at $T=3$. Setting $W=G'$, its linear trajectory is
\begin{equation}\label{eq:traject}
\begin{aligned}
h_1(t,x)&=W(t-x)-W(t+x-1),\\
v_1(t,x)&=W(t-x)+W(t+x-1)\\&-W(t)-W(t-1),\\
s_1(t)&=W(t)+W(t-1),\quad d_1(t)=G(t)+G(t-1).
\end{aligned}
\end{equation}
Direct substitution verifies \eqref{eq:first}, the wall  conditions, and zero initial and final states. The time $T=3$ allows a function on two successive unit intervals and its delayed copy to finish their motion.

Let us define the control subspace $$ \mathcal U:=\{f_G:G\in H^2_0(0,2)\}\subset L^2(0,3). $$ By splitting the time integral into three unit intervals, we get $\|f_G\|_{L^2(0,3)}^2=\|G''\|_{L^2(0,2)}^2+\int_0^1|G''(t)+G''(t+1)|^2 dt.$ Consequently,
\begin{equation}\label{eq:parameter-norm}
\|G''\|_{L^2(0,2)} \le\|f_G\|_{L^2(0,3)} \le\sqrt3\|G''\|_{L^2(0,2)}.
\end{equation}
Invoking the Poincar\'e inequality on $H^2_0(0,2)$, this shows that $\mathcal U$ is closed and that $f_G\mapsto G$ is continuous. Lemma \ref{lem:period} identifies $\Q:\mathcal U\to L^2_0(0,1)$ as a continuous quadratic map. Below, $D\Q(f)$ denotes its Fr\'echet derivative on this control subspace.

The main result constructs a nonzero quadratic return and an
explicit bounded control correction that produces any given mixed second-order response at $T=3$.

\begin{theorem}[A quadratic return and an explicit right inverse]
\label{thm:main}
There exist a nonzero $f_*\in\mathcal U\cap C_c^\infty(0,3)$ and a bounded linear operator $\Rq:L^2_0(0,1)\to\mathcal U$, such that
\begin{equation}\label{eq:reference-return}
\E_1[f_*]=0,\qquad \E_2[f_*]=0.
\end{equation}
For every $q\in L^2_0(0,1)$, the correction $w_q=\Rq q$ satisfies
\begin{equation}\label{eq:main-identities}
\begin{aligned}
\E_1[w_q]&=0,\qquad D\Q(f_*)w_q=q,\\
\Q(f_*+\varepsilon w_q)&=\varepsilon q+\varepsilon^2\Q(w_q)
\quad(\varepsilon\in\mathbb R).
\end{aligned}
\end{equation}
Moreover, for every nonnegative integer $r$, there exists a $C_r>0$, depending only on $r$ and the fixed reference construction, such that
\begin{equation}\label{eq:Hs-bound}
\|w_q\|_{H^r(0,3)}\le C_r\|q\|_{H^r(0,1)},\qquad q\in\mathcal H^r(0,1).
\end{equation}
If, in addition, $q$ is smooth and $1$-periodic, then $w_q\in C_c^\infty(0,3)$.
\end{theorem}
Since by \eqref{eq:main-identities},  $D\Q(f_*)\Rq=\operatorname{Id}_{L^2_0(0,1)}$, $\Rq$ is a right inverse of $D\Q(f_*)$: for each target $q$, the control variation $w_q=\Rq q$ makes the coefficient of $\varepsilon$ in
$\Q(f_*+\varepsilon w_q)$ equal to $q$. Both controls have the delayed form \eqref{eq:f_G}. Section \ref{sec:reference} constructs $G$ with reference acceleration $f_*=f_G$, and Section \ref{sec:inverse} gives correction $w_q=f_{K_q}$ explicitly in \eqref{eq:inverse}--\eqref{eq:control-inv}.

\begin{remark}[Mixed response]\label{rem:variations}
We define $\delta z_1[w]:=z_1[w]$ after fixing $f_*$ and a control variation $w$. Linearity of the first-order equations ensures $z_1[f_*+\varepsilon w]=z_1[f_*]+\varepsilon\delta z_1[w].$ Substitution into the quadratic forcing in
\eqref{eq:second} produces terms proportional to
$1$, $\varepsilon$, and $\varepsilon^2$. The terms proportional to $\varepsilon$ are mixed products of the first-order responses to $f_*$ and $w$. Since \eqref{eq:second} is linear in the second-order unknowns, its solution has the exact decomposition
\begin{equation}\label{eq:mixed-response}
z_2[f_*+\varepsilon w]
=z_2[f_*]+\varepsilon\delta z_2[w]
+\varepsilon^2z_2[w],
\end{equation}
where $\delta z_2[w]$ is the response to these mixed products. Note that $z_2[f_*]$ need not vanish before the final time. This identity holds at every time, and $\delta z_2[w]$ is linear in $w$ for fixed $f_*$. For $w\in\mathcal U$, evaluation at $T=3$ yields $\delta z_2[w] (T,\cdot)|_{(0,1)}=D\Q(f_*)w$. Since $\Q(f_*)=0$, $\Q(f_*+\varepsilon w) =\varepsilon D\Q(f_*)w+\varepsilon^2\Q(w)$. The construction in Section \ref{sec:inverse} proves $D\Q(f_*)w_q=q$, implying \eqref{eq:main-identities}. Note for clarity that parameter $\varepsilon$ varies the reference input, whereas $\alpha$ in \eqref{eq:expansion} scales the physical input.
\end{remark}

\section{Explicit construction}\label{sec:curve}
The proof contains the following steps: 1) reduce the quadratic final-state map to a scalar expression; 2) construct a curve for which this expression is constant; and 3) invert a scalar second-order equation by tangential variation.

\subsection{A scalar formula for the quadratic final state $\Q(\cdot)$}  \label{sec:identity}

Let $\sigma(\cdot)\in L^\infty(\mathbb R)$ be a $2$-periodic function, equal to $1$ on $(0,1)$ and $-1$ on $(-1,0)$. The first and second-order equations \eqref{eq:first}-\eqref{eq:second} become in the unfolded variables \eqref{eq:unfold} 
\begin{equation}\label{eq:unfolded}
 z_{1,t}+z_{1,x}=-\sigma f, \quad z_{2,t}+z_{2,x}=-\partial_x\mathfrak f(z_1,\reflec z_1),
\end{equation}
with zero initial conditions, reflection $(\reflec z)(x):=z(-x)$ and 
\begin{align}
\mathfrak f(a,b):=\frac{3a^2-2ab-b^2}{8}. \label{eq:frak f}
\end{align}
The quadratic flux $\mathfrak f$ is taken from \cite[Lemma 3.1]{CKN}. To see its connection with
the fluid variables, notice that, for $0<x<1$, $\mathfrak f(z_1,\reflec z_1) =h_1v_1+\frac12v_1^2.$ We now integrate the second equation in \eqref{eq:unfolded} along the characteristic $x=X+t$, where $X$ is its position at time $0$. Since $z_2(0,\cdot)=0$, we obtain
\begin{equation}\label{eq:characteristic}
\begin{aligned}
\Q(f)(X+3)&=-\partial_X\Phi_f(X);\\ \Phi_f(X)
&:=\int_0^3 \mathfrak f\bigl(z_1(t,X+t),z_1(t,-X-t)\bigr) dt.
\end{aligned}
\end{equation}
The following identity specializes this characteristic calculation to the delayed controls \eqref{eq:f_G}. 
\begin{lemma}[Scalar reduction]\label{lem:period}
Let us define on $(0,1)$ 
\begin{equation}\label{eq:curve-functional}
\begin{gathered}
g(\xi)=\binom{G(\xi)}{G(\xi+1)},\qquad J=\begin{pmatrix}0&1\\-1&0\end{pmatrix},\\
P_G(\xi)=|g'(\xi)|^2-g'(\xi)^\top Jg(\xi),
\end{gathered}
\end{equation}
for $G\in H^2_0(0,2)$, where $g$ is an auxiliary planar curve parametrized by $\xi \in [0,1]$ and $P_G$ a scalar quadratic expression. Then, $P_G\in H^1(0,1)$ with $P_G(0)=P_G(1)$, and its $1$-periodic extension satisfies
\begin{equation}\label{eq:period}
\Q(f_G)(y)=\tfrac32P_G'(1-y)
\end{equation}
in $L^2_0(0,1)$. It turns out that since $f_G$ is a linear return, it is a quadratic return if and only if $P_G$ is constant.
\end{lemma}
\begin{proof}
Considering $G\in C_c^\infty(0,2)$ and recalling $W=G'$, the first equation in \eqref{eq:unfolded} yields characteristic formula $z_1(t,x)=-\int_{t-1}^t\sigma(x-t+s)G''(s) ds$. The linear trajectory \eqref{eq:traject} vanishes outside $[0,3]$, so the time integral defining $\Phi_{f_G}$ may be extended to $\mathbb R$. We then fix $X$ and set $\tau_j=j-X$, $j\in\mathbb Z$. For $\tau_j<t<\tau_j+1$, the sign changes in the preceding integral occur at $s=\tau_j$ for $x=X+t$ and at $s=2t-\tau_j-1$ for $x=-X-t$ and both points lie in $(t-1,t)$. Splitting at these points, we readily obtain $$ \begin{aligned}
z_1(t,X+t)&=(-1)^j(2A_j-U(t)),\\
z_1(t,-X-t)&=(-1)^j(U(t)-2B_j(t)),
\end{aligned} $$ where $A_j=W(\tau_j)$, $B_j(t)=W(2t-\tau_j-1)$, and $U(t)=W(t)+W(t-1)$. By virtue of \eqref{eq:frak f}, we get $\mathfrak f(2A_j-U,U-2B_j)= \frac{3}{2}A_j^2+A_jB_j-\frac12B_j^2-2A_jU+\frac12U^2.$ Setting $I_j:=\int_{\tau_j-1}^{\tau_j+1}W(s) ds =G(\tau_j+1)-G(\tau_j-1)$, changes of variable yield $$\int_{\tau_j}^{\tau_j+1}A_jB_j dt=\tfrac12A_jI_j, \qquad \int_{\tau_j}^{\tau_j+1}A_jU dt=A_jI_j.$$ The terms containing $A_j$, therefore, sum to $\tfrac32\sum_j(A_j^2-A_jI_j)=\tfrac32P_G(-X)$, with $P_G$ extended periodically. Indeed, let us write $-X=m+\xi$, with $m\in\mathbb Z$ and $0<\xi<1$. Since $G$ and $W$ vanish outside $(0,2)$, only the terms with $\tau_j=\xi$ and $\tau_j=\xi+1$ remain, implying \eqref{eq:curve-functional}.
The remaining values follow by continuity. Furthermore, it is easy to see that 
$$\sum_j\int_{\tau_j}^{\tau_j+1}B_j(s)^2 ds
=\tfrac12\sum_j\int_{\tau_j-1}^{\tau_j+1}W(s)^2 ds
=\int_{\mathbb R}W(s)^2 ds,$$
since almost every point belongs to exactly two of these intervals. This quantity and $\int_{\mathbb R}U^2(s) ds$ are independent of $X$. Hence, $\Phi_{f_G}(X)$ differs from $\tfrac32P_G(-X)$ by a constant. Invoking \eqref{eq:characteristic}, with $y=X+3$ and since $P_G$ is extended with period $1$, we have $\Q(f_G)(y)=\tfrac{3}{2} P_G'(3-y)=\tfrac32P_G'(1-y)$.

If more general $G\in H^2_0(0,2)$ is desired, one may invoke density of $C_c^\infty(0,2)$ in $H^2_0(0,2)$. By the trace conditions, $P_G(0)=P_G(1)=G'(1)^2$. Moreover,
\begin{equation}\label{eq:prod-contin}
P_G'=(2g'-Jg)^\top g''.
\end{equation}
The embedding $H^1(0,1)\hookrightarrow L^\infty(0,1)$ shows that $g\mapsto P_G'$ is continuous from $H^2(0,1;\mathbb R^2)$ to $L^2(0,1)$. The first-order trajectories converge in $C([0,3];H^1(0,1))$, so the spatial derivatives of their quadratic fluxes converge in $C([0,3];L^2(0,1))$. Thus, the characteristic formula passes to the limit in $L^2$ and agrees with \eqref{eq:period}. In particular, \eqref{eq:parameter-norm} and \eqref{eq:prod-contin} imply
$$\begin{aligned}
&\|\Q(f)-\Q(h)\|_{L^2(0,1)}\\
&\quad\le C\bigl(\|f\|_{L^2(0,3)}+\|h\|_{L^2(0,3)}\bigr)\\
&\qquad\quad\times\|f-h\|_{L^2(0,3)},\qquad f,h\in\mathcal U.
\end{aligned}$$
Expanding, next, \eqref{eq:prod-contin} shows that $\mathcal B(f,h) :=\tfrac12\bigl(\Q(f+h)-\Q(f)-\Q(h)\bigr)$ is symmetric and bilinear on $\mathcal U$, with $\|\mathcal B(f,h)\|_{L^2(0,1)} \le C\|f\|_{L^2(0,3)}\|h\|_{L^2(0,3)}.$ It follows $\Q(f+h)=\Q(f)+2\mathcal B(f,h)+\Q(h).$ Since $\|\Q(h)\|_{L^2(0,1)}\le C\|h\|_{L^2(0,3)}^2$, the map $\Q$ is Fr\'echet differentiable on $\mathcal U$, with $$D\Q(f)h=2\mathcal B(f,h) =\Q(f+h)-\Q(f)-\Q(h).$$
\end{proof}

\subsection{Construction of the reference acceleration $f_*$}\label{sec:reference}

The curve $g=(G(\xi),G(\xi+1))^\top$ introduced simply gathers the two successive pieces of the time function $G$. Lemma \ref{lem:period} suggests constructing these pieces together so that
\begin{equation}\label{eq:curve_relation}
|g'|^2-g'^\top Jg=\kappa_*>0.
\end{equation}
The two components of the curve $g$ must assemble into a function $G\in C_c^\infty(0,2)$, while the scalar expression $P_G$ in \eqref{eq:curve-functional} must be constant. We first choose a curve that ensures the required joining, then adjust its parametrization. Choosing $\kappa_*>0$ will also guarantee positivity of the coefficient $a$ in \eqref{eq:a} below.

Let us choose a smooth nondecreasing function $\chi:[0,1]\to[0,1]$ with $\chi=0$ on $[0,1/4]$ and $\chi=1$ on $[1/2,1]$, and set
\begin{equation}\label{eq:gamma}
\gamma(\eta)=
\begin{pmatrix}
(1-\eta)\chi(\eta)\\
-\eta\chi(1-\eta)
\end{pmatrix},\qquad 0\le\eta\le1.
\end{equation}
In particular, $\gamma(\eta)=(0,-\eta)^\top$ on $[0,1/4]$ and $\gamma(\eta)=(1-\eta,0)^\top$ on $[3/4,1]$. These linear expressions will make the two pieces of $G$ match at $t=1$. We then define
\begin{equation}\label{eq:A,B}
A(\eta)=|\gamma'(\eta)|^2,\qquad B(\eta)=\gamma'(\eta)^\top J\gamma(\eta).
\end{equation} We notice that $A\ge1$, since for every $\eta\in[0,1]$, at least one component of
$\gamma'(\eta)$ equals $-1$. Moreover, both $A$ and $B$ are symmetric about $1/2$.

We next set $\Xi(\eta)=\int_0^\eta\rho(s) ds$, for a positive function $\rho$ to be determined, and define the curve $g$ by $g(\Xi(\eta))=\gamma(\eta)$. By the chain rule, $g'(\Xi(\eta))=\gamma'(\eta)/\rho(\eta)$.  Substituting this expression for $g'$ into
\eqref{eq:curve_relation} implies $A/\rho^2-B/\rho=\kappa$. For each $\kappa>0$, its positive solution is expressed as
\begin{equation}\label{eq:rho}
\rho_\kappa(\eta)=
\frac{\sqrt{B(\eta)^2+4\kappa A(\eta)}-B(\eta)}{2\kappa}.
\end{equation}
The choice of $\kappa_*$ comes from the normalization
\begin{equation}\label{eq:normalization}
F(\kappa_*):=\int_0^1\rho_{\kappa_*}(\eta) d\eta=1,
\end{equation}
which ensures that $\Xi$ maps $[0,1]$ onto itself. For completeness, $F(\kappa)=\int_0^1\rho_\kappa(\eta) d \eta$ is continuous and strictly decreasing, since $\partial_\kappa\rho_\kappa =-\frac{\rho_\kappa^2}{\sqrt{B^2+4\kappa A}}<0.$ On $[0,1/4]$, $A=1$ and $B=0$, yielding $F(\kappa)\ge1/(4\sqrt\kappa)\to\infty$ as $\kappa\downarrow0$. Moreover, the bound $\rho_\kappa\le\sqrt{A/\kappa}+|B|/\kappa$ yields $F(\kappa)\to0$ as $\kappa\to\infty$. As a consequence, \eqref{eq:normalization} has a unique positive solution, which can be computed by the bisection method. Substituting $\rho=\rho_{\kappa_*}$, $\Xi(\cdot)$ has a smooth inverse and $g=\gamma\circ\Xi^{-1}$ satisfies \eqref{eq:curve_relation}.

Next, let $k=\sqrt{\kappa_*}$ and $\tau_0=1/(4k)$. Since $\rho=1/k$ on the two end intervals above, the resulting curve satisfies
\begin{equation}\label{eq:ends} g(\xi)=
\begin{cases}
(0,-k\xi)^\top,&0\le\xi\le\tau_0,\\
(k(1-\xi),0)^\top,&1-\tau_0\le\xi\le1.
\end{cases}
\end{equation}
If we assemble
\begin{equation}\label{eq:assembly}
G(t)=
\begin{cases}
g_1(t),&0\le t\le1,\\
g_2(t-1),&1<t\le2,\\
0,&\text{otherwise},
\end{cases}
\end{equation}
then, by \eqref{eq:ends}, we see that $G$ vanishes near $0$ and $2$, and both expressions near $t=1$ are equal to $-k(t-1)$, whence all derivatives match there. It follows that $G\in C_c^\infty(0,2)$.

We now readily get, by the linear trajectory \eqref{eq:traject}, $\E_1[f_*]=0$, if we set $f_*=f_G$. In addition, \eqref{eq:curve_relation} in conjunction with Lemma \ref{lem:period} yield $\Q(f_*)=0$, equivalently, $\E_2[f_*]=0$. Finally, we notice that $\int_0^1f_*(t) dt=G'(1)-G'(0)=-k\ne0$, so $f_*$ is nonzero. This proves \eqref{eq:reference-return}.

\begin{figure}[t]
\centering
\includegraphics[width=0.96\columnwidth]{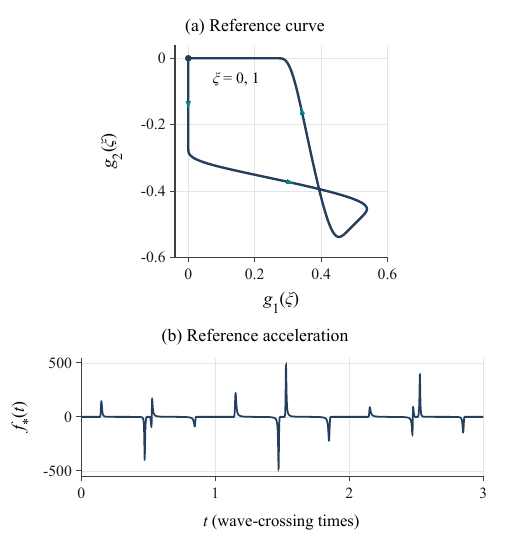}
\caption{Reference construction. The two components of $g$ form the successive pieces of $G$ (arrows indicate increasing curve parameter). The acceleration is $f_*(t)=G''(t)+G''(t-1)$.}
\label{fig:design}
\end{figure}

\subsection{An explicit control correction $w_q$}\label{sec:inverse}

The reference acceleration $f_*$ cancels the first two final-state coefficients. We construct the control correction $w_q$ satisfying \eqref{eq:main-identities} for a smooth $1$-periodic target $q$ with zero mean. We first vary the reference curve $g$ tangentially, as  $g_\varepsilon=g+\varepsilon g'\theta_q$. The tangential variation is useful because the
constant value of $P_G$ eliminates the term
proportional to $\theta_q$, leaving only $\theta_q'$. The scalar function $\theta_q$ will be determined by $q$. One observes that coefficient
\begin{equation}\label{eq:a}
a(\xi)=|g'(\xi)|^2+\kappa_*
\end{equation} is positive and is equal to $2\kappa_*$ near both ends $0$ and $1$, so it extends smoothly with period $1$. The coefficient of $\varepsilon$ in $|g_\varepsilon'|^2-g_\varepsilon'^\top Jg_\varepsilon$ can be written as $2g'^\top(g'\theta_q)'-(g'\theta_q)'^\top Jg-g'^\top J(g'\theta_q)$ $=P_G'\theta_q+(2|g'|^2-g'^\top Jg)\theta_q'=a\theta_q',$
by virtue of $g'^\top Jg'=0$ and $P_G=\kappa_*$. The curve variation $g'\theta_q$ will determine a variation $K_q$ of $G$, and, hence, the control correction $w_q=f_{K_q}$. By virtue of Lemma \ref{lem:period}, the corresponding mixed final response is given by  $\delta z_2[w_q](T,y)=\tfrac32(a\theta_q')'(1-y)$. Requiring this response to equal the target $q(y)$ implies
\begin{equation}\label{eq:theta}
(a\theta_q')'(\xi)=\tfrac23q(1-\xi).
\end{equation}
Performing two integrations, we easily get the following explicit formula for $\theta_q:$
\begin{equation}\label{eq:inverse}
\begin{aligned}
\theta_q(\xi)&:=\int_0^\xi \frac{F_q(\eta)+c_q}{a(\eta)} d\eta,\\
F_q(\xi)&:=\frac23\int_0^\xi q(1-\eta) d\eta, \quad
c_q:=-\frac{\displaystyle\int_0^1 \frac{F_q(\eta)}{a(\eta)} d\eta}{\displaystyle\int_0^1\frac{1}{a(\eta)} d\eta}.
\end{aligned}
\end{equation} 
Since $q$ has zero mean, we have $F_q(0)=F_q(1)=0$, and by the choice of $c_q$, we have $\theta_q(0)=\theta_q(1)=0$, and $\theta_q'(0)=\theta_q'(1)=c_q/(2\kappa_*)$. In addition, all higher derivatives match because $q$ and $a$ are smooth and $1$-periodic. Thus, $\theta_q$ above is smooth and $1$-periodic.

Let us now assemble the two components of $g'\theta_q$ into
\begin{equation}\label{eq:correction-assmbl}
K_q(t)=
\begin{cases}
g_1'(t)\theta_q(t),&0\le t\le1,\\
g_2'(t-1)\theta_q(t-1),&1<t\le2,\\
0,&\text{otherwise},
\end{cases}
\end{equation}
and define the control correction by
\begin{equation}\label{eq:control-inv}
w_q(t)=K_q''(t)+K_q''(t-1),\qquad \Rq q=w_q.
\end{equation}
By \eqref{eq:ends}, $K_q$ vanishes near $0$ and $2$.
Near $t=1$, its two expressions are $-k\theta_q(t)$ and $-k\theta_q(t-1)$, so periodicity again makes all derivatives match. These guarantee $K_q\in C_c^\infty(0,2)$ and $\E_1[w_q]=0$. Moreover, the scalar quadratic expression \eqref{eq:curve-functional} satisfies the exact identity
\begin{equation}\label{eq:variation}
P_{G+\varepsilon K_q}=\kappa_*+\varepsilon a\theta_q' +\varepsilon^2P_{K_q}.
\end{equation}
Since $f_*+\varepsilon w_q=f_{G+\varepsilon K_q}$,
by Lemma \ref{lem:period} and \eqref{eq:variation} 
$\Q(f_*+\varepsilon w_q)(y)=\tfrac32P_{G+\varepsilon K_q}'(1-y)$ $=\tfrac{3\varepsilon}{2}(a\theta_q')'(1-y)+\varepsilon^2\Q(w_q)(y)$ $=\varepsilon q(y)+\varepsilon^2\Q(w_q)(y),$ where the last equality follows from \eqref{eq:theta}. The coefficient of $\varepsilon$ is, therefore,
$D\Q(f_*)w_q=q$. Combined with $\E_1[w_q]=0$
proved above, this establishes \eqref{eq:main-identities}. The construction is linear in $q$.

The same integral formulas apply to $q\in L^2_0(0,1)$. Since $a(\cdot)$ is smooth and positive and $|c_q|\le\|F_q\|_{L^\infty(0,1)}$, they ensure bound $$ \|\theta_q\|_{H^{r+2}(0,1)} \le C_r\|q\|_{H^r(0,1)},\qquad q\in\mathcal H^r(0,1),$$ for some constant $C_r>0$, and for every nonnegative integer $r.$ The matching traces of $\theta_q$ through derivative order $r+1$, combined with \eqref{eq:ends}, ensure that the two pieces of $K_q$ join in $H^{r+2}$ at $t=1$. Since $K_q$ vanishes near $0$ and $2$, its zero extension belongs to $H^{r+2}(\mathbb R)$. For $r=0$, matching values and first derivatives ensures that $K_q''$ is an $L^2$ function, without point impulses. Therefore, there exists $\widetilde C_r>0,$ such that $$\|K_q\|_{H^{r+2}(\mathbb R)}+\|w_q\|_{H^r(0,3)} \le\widetilde C_r\|q\|_{H^r(0,1)}.$$ Note that constants $C_r, \widetilde C_r$ depend only on $r$ and on the fixed curve $g$, and not on $q$. We then see that approximation by smooth $1$-periodic functions with zero mean, in conjunction with the continuity proven in Lemma \ref{lem:period}, extends \eqref{eq:main-identities} to $L^2_0(0,1)$. By this, the proof of Theorem \ref{thm:main} is complete.

\subsection{Exact correction within the quadratic equations}
The final step to reveal the exact correction is to invoke a contraction argument.  In fact, the quadratic remainder $\Q(w_q)$ in \eqref{eq:main-identities} can be removed for sufficiently small targets $q$ by a further correction ($\varepsilon=1$).
\begin{corollary}\label{cor:quadratic-local}
There exist $\delta,C>0$, such that every target $q\in L^2_0(0,1)$ with $\|q\|_{L^2(0,1)}<\delta$ admits an exact correction for the quadratic equations $\widehat w_q\in\mathcal U$ satisfying
\begin{equation}\label{eq:quadratic-local}
\begin{aligned}
\Q(f_*+\widehat w_q)&=q,\\
\|\widehat w_q\|_{L^2(0,3)}&\le C\|q\|_{L^2(0,1)}.
\end{aligned}
\end{equation}
The corrected control $f_*+\widehat w_q $ remains a linear return.
\end{corollary}
\begin{proof}
Let $p\in L^2_0(0,1)$ be an auxiliary function and set $N(p)=\Q(\Rq p)$. We seek $p$ such that $p+N(p)=q$. Since map $N$ is continuous and quadratic, there exists an $M>0$, such that $\|N(p_1)-N(p_2)\|_{L^2(0,1)}$ $\le M\bigl(\|p_1\|_{L^2(0,1)}+\|p_2\|_{L^2(0,1)}\bigr) \|p_1-p_2\|_{L^2(0,1)}. $ Taking $0<\|q\|_{L^2(0,1)}<1/(8M)$, the map $p\mapsto q-N(p)$ maps the closed ball of radius $2\|q\|_{L^2(0,1)}$ into itself
and has a Lipschitz constant less than $1/2$ there.
We define $$p_0=q, \qquad p_{n+1}=q-N(p_n), \quad n\ge0.$$ The first control approximation is $f_*+w_q$, since $\Rq p_0=w_q$. Banach fixed-point theorem ensures convergence in $L^2_0(0,1)$ to a fixed point $p_\infty$ and we set $$\widehat w_q=\Rq p_\infty.$$ Then, \eqref{eq:main-identities}, with $\varepsilon=1$, yields $\Q(f_*+\widehat w_q)=q$. Each iterate provides a control whose error satisfies
$\Q(f_*+\Rq p_n)-q=p_n-p_{n+1}.$ Since $\|p_1-p_0\|_{L^2(0,1)} \le M\|q\|_{L^2(0,1)}^2$, contraction ensures
$ \|\Q(f_*+\Rq p_n)-q\|_{L^2(0,1)}
\le 2^{-n}M\|q\|_{L^2(0,1)}^2.$
We notice that every evaluation of $N$ uses the explicit formula $$N(p)(y)=\tfrac32P_{K_p}'(1-y)$$ from Lemma \ref{lem:period}, with $K_p$ constructed by \eqref{eq:inverse}--\eqref{eq:correction-assmbl} with $q$ replaced by $p$. Finally, by $\|p_\infty\|_{L^2(0,1)}\le2\|q\|_{L^2(0,1)}$, one gets $\|\widehat w_q\|_{L^2(0,3)} \le 2 C_0\|q\|_{L^2(0,1)}$, where $C_0$ is the constant in \eqref{eq:Hs-bound} for $r=0$. For $q=0$, one may take $\widehat w_q=0$.
\end{proof}
The correction $w_q=\Rq q$ is given explicitly and
produces the mixed final response $q$. The exact correction $\widehat w_q$ is obtained as the limit in $L^2(0,3)$ of the preceding iteration and satisfies $\Q(f_*+\widehat w_q)=q$. Thus, the quadratic final-state equation is solved constructively.

\begin{remark}[Scope of the quadratic result]
A quadratic return need not be an exact return of
the nonlinear system. Corollary \ref{cor:quadratic-local} solves only the quadratic final-state equation. Nonlinear controllability at $T=3$ requires additional estimates for the full solution map and its linearization, which are not supplied by \eqref{eq:Hs-bound}.
\end{remark}

\section{A Sylvester equation for the missing modes}
\label{sec:sylvester}
Section \ref{sec:inverse} constructs the control correction $w_q$ through the auxiliary curve $g$. We now express the same mixed final response using a Sylvester transformation. This transformation removes the coupling between
$\delta z_1$ and $\delta z_2$ and expresses its effect through the input $w$. It also identifies the spatial derivative required by the change of state coordinates. For systems in cascade, indirect controllability of conservative PDEs with a periodic evolution group was studied in \cite{RosierTeresa2011}, while backstepping combined with high-gain observation was used for output-feedback stabilization of linear Korteweg--de Vries equations in \cite{KitsosCerpa2021}. Related differential transformations are developed in \cite[Sec.~3.3.2]{KitsosThesis}, \cite{KBP,KitsosBajodekBaudouin}, while a modal approach is used in \cite{KitsosFridman2024}. 

\subsection{The variation cascade}
We use the variations defined in Remark \ref{rem:variations}, writing $\delta z_i:=\delta z_i[w]$, $i=1,2$ and $z_*:=z_1[f_*]$, for control variation $w\in L^2(0,3)$. For the control correction $w=w_q$ of Theorem \ref{thm:main}, we have $\delta z_1(T,\cdot)=0$ and $\delta z_2(T,\cdot)|_{(0,1)}=q$.  All spatial functions below are extended with period $2$ from $(-1,1)$. We define the translation and multiplication operators by $(V(t)z)(x)=z(x-t), M_cz=cz,$ recalling the reflection $\reflec$ defined in Section \ref{sec:identity}. The transformation used below acts from the Fourier modes with odd indices of $\delta z_1$ to the  modes with even indices of $\delta z_2$. We, therefore, distinguish the following two subspaces:
\begin{align*}
H_o^0&=\{z\in L^2(-1,1):z(x+1)=-z(x)\},\\
H_e^0&=\left\{z\in L^2(-1,1):z(x+1)=z(x),\
\int_{-1}^1z(x) dx=0\right\},
\end{align*}
with the translation identities holding almost everywhere. Let $\widetilde H^1 :=\{z\in H^1(-1,1):z(-1)=z(1)\}$ and  $H_j^1=H_j^0\cap\widetilde H^1$, for $j=o,e$. For clarity, restriction to $(0,1)$ identifies $H_e^r$ with $\mathcal H^r(0,1)$ for $r=0,1$, with norms differing by a factor $\sqrt2$. We consider also the transport generators  $A_j=-\partial_x$, with domain $H_j^1$ in $H_j^0$; $j=o,e$, and the  input operator $\mathsf B w=-\sigma w$, where $\sigma$ is introduced in Section \ref{sec:identity}.

Let $(h_*,v_*)$ be the first-order response to the reference acceleration $f_*$ of Section \ref{sec:reference}, whose unfolding is $z_*=z_1[f_*]$, and set
\begin{equation}\label{eq:sylv-ab}
a_*=\frac{3z_*-\reflec z_*}{4},\qquad b_*=-\frac{z_*+\reflec z_*}{4}.
\end{equation}
These coefficients are Lipschitz continuous on $[0,3]\times[-1,1]$, are antiperiodic under translation by $1$, and vanish at $t=0$. On $(0,1)$, $a_*=v_*+h_*/2$ and $b_*=-h_*/2$. Expanding the quadratic flux \eqref{eq:frak f}, its coefficient of $\varepsilon$ is $a_*\delta z_1+b_*\reflec\delta z_1$. Therefore, \eqref{eq:unfolded} yields the cascade system
\begin{equation}\label{eq:sylv-cascade}
\begin{aligned}
\partial_t\delta z_1&=A_o\delta z_1+\mathsf B w,\\
\partial_t\delta z_2&=A_e\delta z_2+C(t)\delta z_1,
\end{aligned}
\end{equation}
with zero initial conditions and coupling operator $C(t)=-\partial_x(M_{a_*}+M_{b_*}\reflec)$, noting that $C(t):H_o^1\to H_e^0$ is bounded. Since reference acceleration $f_*$ is fixed, coefficients  $a_*,b_*$ are known, and \eqref{eq:sylv-cascade} is linear. The transformation $y_2=\delta z_2+S(t)\delta z_1$ removes the internal coupling when time-dependent Sylvester equation
\begin{equation}\label{eq:sylvequation}
\dot S+SA_o-A_eS=-C,\qquad S(0)=0
\end{equation}
holds. This identity is understood in the sense of distributions. 

\subsection{An explicit differential transformation}

The coefficient $\psi$ depends only on $f_*$. It is obtained by integrating $a_*$ along characteristics and determines the differential part of $S(t)$.
\begin{proposition}\label{prop:sylvester}
Let us define 
\begin{align}
\psi(t,x)&=\int_0^t a_*(\tau,x-t+\tau) d\tau,
\label{eq:sylv-theta}\\
\ell&=\tfrac12(\partial_xb_*-\partial_tb_*),\nonumber\\
\mathsf K(t)&=\int_0^tV(t-\tau)M_{\ell(\tau)} \reflec V(\tau-t) d\tau, \label{eq:sylv-K}\\
S(t)&=\partial_xM_{\psi(t)} +\tfrac12M_{b_*(t)}\reflec+\mathsf K(t).
\label{eq:sylv-S}
\end{align}
Then, $S(t):H_o^1\to H_e^0$ is uniformly bounded and strongly continuous for $0\le t\le3$, and satisfies \eqref{eq:sylvequation} in the sense of distributions. Whenever $\psi(t,\cdot)\ne0$, this operator has no bounded extension $H_o^0\to H_e^0$.

Assuming $w\in L^2(0,3)$, the cascade states satisfy $\delta z_1\in C([0,3];H_o^1)$ and $\delta z_2\in C([0,3];H_e^0)$. Moreover,
\begin{equation}\label{eq:sylv-target}
\dot y_2=A_ey_2+S(t)\mathsf B w,\qquad y_2=\delta z_2+S(t)\delta z_1
\end{equation}
holds in periodic spatial distributions. Above, $S(t)\mathsf B$ uses the distributional extension of \eqref{eq:sylv-S} (its Duhamel integral has a representative in $C([0,3];H_e^0)$).
\end{proposition}
\begin{proof}
We set for convenience $D=\partial_x$ and $[D,S]=DS-SD$. By definition, $\psi$ satisfies a transport equation with source $a_*$ and zero initial value, namely, $\psi_t+\psi_x=a_*$. This, combined with $\reflec D=-D\reflec$, yields
\begin{align*}
\partial_t(DM_\psi)+[D,DM_\psi]&=DM_{a_*},\\
\partial_t(M_{b_*}\reflec/2)&+[D,M_{b_*}\reflec/2]\\
 &=\tfrac12M_{\partial_tb_*+\partial_xb_*}\reflec+M_{b_*}D\reflec,\\ \dot{\mathsf K}+[D,\mathsf K]&=M_\ell\reflec.
\end{align*}
Their sum is then equal to $D(M_{a_*}+M_{b_*}\reflec)=-C$. The definitions \eqref{eq:sylv-theta} and \eqref{eq:sylv-K} induce $\psi(0,\cdot)=0$ and $\mathsf K(0)=0$. Since $z_*(0,\cdot)=0$, \eqref{eq:sylv-ab} also gives $b_*(0,\cdot)=0$. Hence, \eqref{eq:sylv-S} yields $S(0)=0$. Note that derivatives are understood by testing against smooth $2$-periodic functions. Equivalently, integrating along the transport group yields \begin{equation}\label{eq:sylv-integral}
S(t)=-\int_0^tV(t-\tau)C(\tau)V(\tau-t) d\tau. \end{equation} For $z\in H_o^1$, each factor in the products $a_*z$ and $b_*\reflec z$ changes sign under translation by $1$. These products are, therefore, $1$-periodic, and their spatial derivatives have zero mean. Thus, $C(t)z\in H_e^0$. Since translations preserve the mode subspaces, \eqref{eq:sylv-integral} maps $H_o^1$ into $H_e^0$. The Lipschitz coefficients establish a uniform bound. Factoring $V(t)$ and $V(-t)$ outside this integral also proves strong continuity of $S(t):H_o^1\to H_e^0$, since the remaining integral is Lipschitz in operator norm.

The last two terms of \eqref{eq:sylv-S} are bounded on $L^2(-1,1)$. Using the complexification of $H_o^0$ we write, for odd $n$, $z_n(x)=2^{-1/2}e^{in\pi x}$.
Then, $\|z_n\|_{L^2(-1,1)}=1$ and $$ \|\partial_x(\psi z_n)\|_{L^2(-1,1)} \ge\frac{|n|\pi\|\psi\|_{L^2(-1,1)}-\|\psi_x\|_{L^2(-1,1)}}{\sqrt2}.$$ If $\psi\ne0$, this grows with $|n|$ and cannot be canceled by the bounded terms. Hence, $S(t)$ has no bounded $L^2$ extension.

The first equation in \eqref{eq:sylv-cascade} yields $\delta z_1(t,x)=-\int_0^t\sigma(x-t+\tau)w(\tau) d\tau,$ showing that its spatial derivative is a finite sum of translated restrictions of $w$, because $\sigma$ is piecewise constant. Continuity in time follows from continuity of $t\mapsto w\mathbf1_{(0,t)}$ in $L^2(\mathbb R)$, after zero extension, and strong continuity of translations. The second equation in \eqref{eq:sylv-cascade} then yields, for some $C>0,$ $\|\delta z_1\|_{C([0,3];H^1(-1,1))}+\|\delta z_2\|_{C([0,3];L^2(-1,1))} \le C\|w\|_{L^2(0,3)}.$ Taking $w$ smooth, if we differentiate $y_2=\delta z_2+S\delta z_1$ in the sense of distributions, \eqref{eq:sylv-target} and $$y_2(t)=\int_0^tV(t-\tau)S(\tau)\mathsf B w(\tau) d\tau $$ are proven. This integral is initially understood as a distribution. Formula \eqref{eq:sylv-S} also defines $S(t)z$ as a  distribution for $z\in L^2(-1,1)$. Performing integration by parts, we get $|\langle S(t)z,\varphi\rangle| \le C\|z\|_{L^2(-1,1)}\|\varphi\|_{H^1(-1,1)}$ for smooth periodic test functions $\varphi$. Approximating any $w\in L^2(0,3)$ by smooth controls, the integral converges in distributions while $\delta z_2+S\delta z_1$ converges in $C([0,3];L^2(-1,1))$. The limits agree, proving the asserted regularity of the input integral.
\end{proof}

\subsection{The spatial derivative and the control correction $w_q$}

The first term of \eqref{eq:sylv-S} acts as $\partial_x(\psi\delta z_1) =\psi\partial_x\delta z_1+\psi_x\delta z_1.$ Thus, this change of coordinates requires one spatial derivative of $\delta z_1$ to produce an $L^2$ state. This requirement comes from the spatial derivative in $C(t)$, already present in \eqref{eq:unfolded} without requiring an additional time derivative
of the scalar control $w$. The reflected part behaves differently, since reflection reverses the transport direction, and $[D,M_{b_*}\reflec]=M_{\partial_xb_*}\reflec+2M_{b_*}D\reflec$. The factor $2$ explains the coefficient $1/2$ in \eqref{eq:sylv-S}, which leaves only bounded terms for this part of the coupling.

For the constructed reference, $\psi$ is not identically zero. Indeed, using \eqref{eq:traject} and the unfolding \eqref{eq:unfold}, we obtain
\begin{equation}\label{eq:sylv-wall-coefficient} \psi(t,0)=\tfrac34G(t),\qquad 0\le t\le1.
\end{equation}
Since $G$ is not identically zero on $(0,1)$, the differential term cannot be omitted from this transformation.

The triangular map $(\delta z_1,\delta z_2)\mapsto (\delta z_1,\delta z_2+S\delta z_1)$ is a bounded isomorphism on $H_o^1\times H_e^0$. The proposition shows why it cannot instead be used as a bounded map on $H_o^0\times H_e^0$ when $\psi\ne0$. This constitutes the derivative requirement of the state transformation, as in the differential constructions of \cite{KitsosThesis,KBP}, which is compatible with $w\in L^2(0,3)$, since its first-order response $\delta z_1$ belongs to $H_o^1$.

Taking $w=w_q=\Rq q$, the linear return implies $\delta z_1(3,\cdot)=0$. Hence, the transformed and original mixed responses coincide at the final time. Restricting to $(0,1)$ yields
\begin{equation}\label{eq:sylv-inverse-link}
D\Q(f_*)w_q=\left.\int_0^3 V(3-t)S(t)\mathsf B w_q(t) dt\right|_{(0,1)}=q.
\end{equation} 
In this way, the Sylvester equation identifies a transformed input $S(t)\mathsf B w(t)$ acting on the missing modes, and the curve calculation provides a right inverse for its final-state map.

In the preceding construction, the same differentiated quadratic forcing was integrated along characteristics in \eqref{eq:characteristic}. For the tangential correction, its final effect is $ \delta z_2(3,y)=\tfrac32(a\theta_q')'(1-y),$ by \eqref{eq:period} and \eqref{eq:variation}, with $w=w_q$. Note here that $a(\xi)$ is the scalar coefficient in \eqref{eq:a}, distinct from $a_*(t,x)$ in \eqref{eq:sylv-ab}. Solving \eqref{eq:theta} ensures  $\theta_q\in H^{r+2}(0,1)$ for $q\in\mathcal H^r(0,1)$. The assembly in \eqref{eq:correction-assmbl} implies $K_q\in H^{r+2}(\mathbb R)$, and the two time derivatives in \eqref{eq:control-inv} ensure $w_q\in H^r(0,3)$. Note that differential elimination and derivative loss also arise in controllability studies for hyperbolic systems with a reduced number of controls \cite{Alabau-BoussouiraCoronOliv2017}. In the present case, the derivative requirement concerns the state transformation $S$, while the final-state correction $q\mapsto w_q$ preserves regularity in the Sobolev space chosen here.

\section{Numerical illustration}\label{sec:numerics}
The construction of Section \ref{sec:reference} is implemented with $$\chi(\eta)=\frac{b_0(4\eta-1)}{b_0(4\eta-1)+b_0(2-4\eta)},$$ where $b_0(z)=e^{-1/z}$ for $z>0$ and $b_0(z)=0$ otherwise. Applying the bisection method in \eqref{eq:normalization}, one gets $\kappa_*\simeq3.5995$. We then form $g=\gamma\circ\Xi^{-1}$, assemble $G$ using \eqref{eq:assembly}, and compute $f_*=f_G$ using \eqref{eq:f_G}. Figure \ref{fig:design} shows the reference curve and acceleration.

\begin{figure}[t]
\centering
\includegraphics[width=0.96\columnwidth]{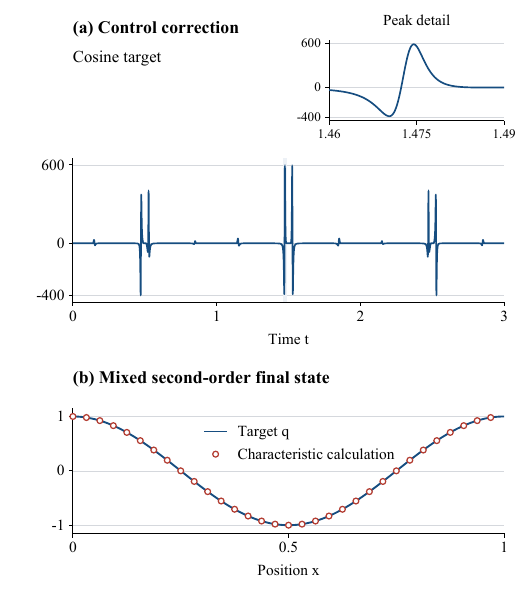}
\caption{Correction for $q(x)=\cos(2\pi x)$. Top: $w_q$ and a detail of one narrow peak. Bottom: the target and mixed final state computed from characteristic integration. Circles mark every fourth spatial sample.}
\label{fig:correction}
\end{figure}

Assuming a desired target $q(x)=\cos(2\pi x)$, \eqref{eq:recovered} yields a mixed height response $\cos(2\pi x)$ and zero mixed velocity. Symmetry of $A$ and $B$ in \eqref{eq:A,B} enforces $\rho(1-\eta)=\rho(\eta)$ and $\Xi(1-\eta)=1-\Xi(\eta)$, hence, $a(1-\xi)=a(\xi)$. Formula \eqref{eq:inverse}, therefore, yields $F_q(\xi)=\frac{\sin(2\pi\xi)}{3\pi}, c_q=0,$
and $$\theta_q(\xi)=\frac{1}{3\pi}\int_0^\xi\frac{\sin(2\pi\eta)}{a(\eta)} d\eta.$$ We assemble $K_q$ and compute $w_q$ using \eqref{eq:correction-assmbl}--\eqref{eq:control-inv}; see Fig. \ref{fig:correction}.

To check the mixed-response identity independently of \eqref{eq:period}, we integrate the mixed flux
$a_*\delta z_1[w_q]+b_*\reflec\delta z_1[w_q]$ along the characteristics in \eqref{eq:characteristic}, using Gauss--Legendre quadrature separately between consecutive wall crossings. The negative spatial derivative of the resulting periodic function is computed by Fourier spectral differentiation on the uniform spatial grid. At $x_j=j/128$, one compares the resulting $\delta z_2^{\mathrm{num}}(3,x_j)$ with $q(x_j)$ using $e_{\mathrm{samp}}=\max_{0\le j<128}|\delta z_2^{\mathrm{num}}(3,x_j)-q(x_j)|.$ Table \ref{tab:refinement} reports the sampled error. Refining both the curve representation and characteristic quadrature reduces this error, keeping the spatial grid fixed. This calculation checks $D\Q(f_*)w_q=q$ (the iterative computation of the exact correction $\widehat w_q$ in Corollary \ref{cor:quadratic-local} is not illustrated here).
\begin{table}[t]
\vspace*{5pt} 
\centering
\caption{Sampled characteristic error ($N_x=128$).}
\label{tab:refinement}
\begin{tabular}{rrr}
\toprule
Curve points & Gauss points/interval & $e_{\mathrm{samp}}$\\
\midrule
131073 & 2048 & $1.38\times10^{-4}$\\
262145 & 4096 & $3.63\times10^{-5}$\\
\bottomrule
\end{tabular}
\end{table}
The acceleration plots show the reference $f_*$
and correction $w_q$. The corresponding input is $u=\alpha(f_*+\varepsilon w_q)$, as obtained in the previous sections.

\section{Conclusion}\label{sec:conclusion}
At $T=3$, the planar construction yields a nonzero quadratic return and an explicit right inverse of $D\Q(f_*)$, preserving regularity in Sobolev spaces on all missing modes. Extension to the full nonlinear model at $T=3$ requires further analysis.

\end{document}